\documentclass{etds}   \usepackage{amsfonts} \usepackage{latexsym} \usepackage{graphicx} \usepackage{varioref}   

\newcommand{\D}{\displaystyle}   %
\newcommand{\breath}{} 

\newtheorem{thm}{Theorem}
\newcounter{claimcount}[thm] 

\newtheorem{prop}[thm]{Proposition}

\newtheorem{lemma}[thm]{Lemma}

\newtheorem{cor}[thm]{Corollary} 

\newtheorem{defn}[thm]{Definition}

\newcommand{\Theorem}[2]{\begin{thm}#1  #2 \end{thm}}
\newcommand{\Proposition}[2]{\begin{prop}#1  #2 \end{prop}}
\newcommand{\Lemma}[2]{\begin{lemma}#1  #2 \end{lemma}}
\newcommand{\Corollary}[2]{\begin{cor}#1  #2 \end{cor}} 
\newcommand{\thmfont}[1]{{\sl #1}}    
\newcommand{\Definition}[2]{\begin{defn}{\em #1:} \\ #2 \end{defn}} 

\newcommand{\proposition}[2]{\Proposition{#1}{#2}}
\newcommand{\Claim}[1]{\refstepcounter{claimcount}                {\sc Claim \theclaimcount: \ }\thmfont{ #1}} 
\newcommand{\claim}{\Claim}
\newcommand{\bprf}[1][Proof:]{\begin{list}{} 			{\setlength{\leftmargin}{1em} 			\setlength{\rightmargin}{0em}}                         \item {\sc \hspace{-1.1em}  #1 \ \ }} 
\newcommand{\eprf}{\end{list}} 
\newcommand{\bprfof}[1]{\begin{list}{} 			{\setlength{\leftmargin}{1em} 			\setlength{\rightmargin}{0em}}                         \item {\sc \hspace{-1em}  Proof of #1: \ \ }} 
\newcommand{\bthmprf}{\bprf}
\newcommand{\bclaimprf}{\bprf}
\newcommand{\ethmprf}{ {\sc \hrulefill $\Box$ }\eprf  \breath } 
\newcommand{\eclaimprf}{ {\sc \dotfill~$\Box$~[Claim~\theclaimcount] }\eprf}
\newcommand{\QED}{\hrulefill\ensuremath{\Box}}
\newcommand{\qed}{\QED}     
\newcommand{\beq}{\begin{eqnarray*}}
\newcommand{\eeq}{\end{eqnarray*}} 
\newcommand{\beqn}{\begin{equation}}
\newcommand{\eeqn}{\end{equation}} 
\newcommand{\blist}{\begin{enumerate}}
\newcommand{\elist}{\end{enumerate}} 
\newcommand{\bitem}{\begin{itemize}}
\newcommand{\eitem}{\end{itemize}} 
\newcommand{\Cesaro}{Ces\`aro }
\newcommand{\done}{\ensuremath{\mathsf{ 1\!\!1}}} 
\newcommand{\dB}{\ensuremath{\mathbb{B}}}
\newcommand{\dC}{\ensuremath{\mathbb{C}}}
\newcommand{\dG}{\ensuremath{\mathbb{G}}}
\newcommand{\dJ}{\ensuremath{\mathbb{J}}}
\newcommand{\dM}{\ensuremath{\mathbb{M}}}
\newcommand{\dN}{\ensuremath{\mathbb{N}}}
\newcommand{\dP}{\ensuremath{\mathbb{P}}}
\newcommand{\dT}{\ensuremath{\mathbb{T}}}
\newcommand{\dU}{\ensuremath{\mathbb{U}}}
\newcommand{\dV}{\ensuremath{\mathbb{V}}}
\newcommand{\dZ}{\ensuremath{\mathbb{Z}}}   %
\newcommand{\bC}{\ensuremath{\mathbf{ C}}}
\newcommand{\bL}{\ensuremath{\mathbf{ L}}}
\newcommand{\ba}{\ensuremath{\mathbf{ a}}}
\newcommand{\bb}{\ensuremath{\mathbf{ b}}}
\newcommand{\bc}{\ensuremath{\mathbf{ c}}}
\newcommand{\bff}{\ensuremath{\mathbf{ f}}}
\newcommand{\bg}{\ensuremath{\mathbf{ g}}}
\newcommand{\Bh}{\ensuremath{\mathbf{ h}}}
\newcommand{\bi}{\ensuremath{\mathbf{ i}}}
\newcommand{\bk}{\ensuremath{\mathbf{ k}}}
\newcommand{\Bm}{\ensuremath{\mathbf{ m}}}
\newcommand{\bq}{\ensuremath{\mathbf{ q}}}
\newcommand{\bs}{\ensuremath{\mathbf{ s}}}
\newcommand{\bw}{\ensuremath{\mathbf{ w}}}
\newcommand{\bx}{\ensuremath{\mathbf{ x}}}
\newcommand{\bgam }{\mbox{{\boldmath $\gamma$}}}
\newcommand{\bsig }{\mbox{{\boldmath $\sigma$}}}
\newcommand{\bxi}{\mbox{{\boldmath $\xi $}}}
\newcommand{\bnu}{\mbox{{\boldmath $\nu $}}} 
\newcommand{\bchi}{\mbox{{\boldmath $\chi $}}}
\newcommand{\sA}{\ensuremath{\mathcal{ A}}}
\newcommand{\sB}{\ensuremath{\mathcal{ B}}}

\newcommand{\sG}{\ensuremath{\mathcal{ G}}}
\newcommand{\sH}{\ensuremath{\mathcal{ H}}}
\newcommand{\sL}{\ensuremath{\mathcal{ L}}}
\newcommand{\sM}{\ensuremath{\mathcal{ M}}}
\newcommand{\sP}{\ensuremath{\mathcal{ P}}}
\newcommand{\sQ}{\ensuremath{\mathcal{ Q}}}
\newcommand{\sS}{\ensuremath{\mathcal{ S}}}
\newcommand{\sU}{\ensuremath{\mathcal{ U}}}
\newcommand{\brbP}{\ensuremath{\overline{\underline{\mathbf{ P}}}}}
\newcommand{\brbQ}{\ensuremath{\overline{\underline{\mathbf{ Q}}}}}
\newcommand{\gF}{\ensuremath{\mathfrak{ F}}}
\newcommand{\gG}{\ensuremath{\mathfrak{ G}}}
\newcommand{\gH}{\ensuremath{\mathfrak{ H}}}
\newcommand{\gf}{\ensuremath{\mathfrak{ f}}}
\newcommand{\gothg}{\ensuremath{\mathfrak{ g}}}
\newcommand{\gh}{\ensuremath{\mathfrak{ h}}}
\newcommand{\alp }{\ensuremath{\alpha}}
\newcommand{\del }{\ensuremath{\delta}}
\newcommand{\eps }{\ensuremath{\epsilon}}
\newcommand{\h}[1]{{\ensuremath{\widehat{#1}}}} 
\newcommand{\hsA}{\ensuremath{{\widehat{\mathcal{ A}}}}}
\newcommand{\hmu }{\ensuremath{\widehat{\mu}}}
\newcommand{\hnu }{\ensuremath{\widehat{\nu}}}
\newcommand{\hphi}{\ensuremath{{\widehat{\phi }}}}
\newcommand{\fe}{\ensuremath{\mathsf{ e}}}
\newcommand{\fk}{\ensuremath{\mathsf{ k}}}
\newcommand{\fm}{\ensuremath{\mathsf{ m}}}
\newcommand{\fn}{\ensuremath{\mathsf{ n}}}
\newcommand{\fu}{\ensuremath{\mathsf{ u}}}
\newcommand{\fv}{\ensuremath{\mathsf{ v}}}
\newcommand{\vN}{\ensuremath{{\vec{N}}}}
\newcommand{\vm}{\ensuremath{{\vec{m}}}}
\newcommand{\vn}{\ensuremath{{\vec{n}}}}
\newcommand{\vq}{\ensuremath{{\vec{q}}}}
\newcommand{\lb}{\ensuremath{\left}}
\newcommand{\rb}{\ensuremath{\right}} 
\newcommand{\maketall}{\rule[-0.5em]{0em}{1em}}       
\newcommand{\Implies}{\mbox{$\Longrightarrow$}}
\newcommand{\map}{\ensuremath{\longrightarrow}}
\newcommand{\goto}{\ensuremath{\rightarrow}}
\newcommand{\into}{\ensuremath{\map}}
\newcommand{\statement}[1]{\lb( \ \maketall       \begin{minipage}{40em}       \begin{tabbing}         #1        \end{tabbing}      \end{minipage} \ \rb)}     
\newcommand{\oo}{\ensuremath{\infty}}        
\newcommand{\X}{\ensuremath{\times}}
\newcommand{\x}{\ensuremath{\X}}
\newcommand{\tensor}{\ensuremath{\otimes}}
\newcommand{\Tensor}{\ensuremath{\bigotimes}}
\newcommand{\dirsum}{\ensuremath{\oplus}}
\newcommand{\Union}{\ensuremath{\bigcup}}
\newcommand{\intsct}{\ensuremath{\cap}}
\newcommand{\set}[2]{\ensuremath{\left\{ #1 \; ; \; #2 \right\} }} 
\newcommand{\seq}[2]{\ensuremath{ \lb\{#1 |_{_{{#2}}} \rb\} }}     
\newcommand{\pbinom}[3]{\ensuremath{  \left[ { #1 } \atop { #2 } \right]_{#3} }} 
\newcommand{\norm}[2]{\ensuremath{\left\| #1 \right\|_{{#2}} }   }
\newcommand{\inn}[1]{\ensuremath{\left\langle #1 \right\rangle }}       
\newcommand{\Id}[1]{\ensuremath{\mathbf{ Id}_{{#1}}}}
\newcommand{\pr}[1]{\ensuremath{\mathbf{ pr}_{{#1}}}}
\newcommand{\chr}[1]{\ensuremath{{\done}_{{#1}}}} 
\newcommand{\shift}[1]{ {{\bsig}^{#1}} }   
\newcommand{\Prob}[1]{{\sf Prob}\lb[ #1 \rb] } 
\newcommand{\goesto}[2]{\ensuremath{ -\!\!\!-\!\!\!-\!\!\!-\!\!\!\!\!\!\!\!\!\!\!  ^{{\scriptscriptstyle #2}}_{{\scriptscriptstyle #1}}   \!\!\!\!\!\!\!\!\!\longrightarrow }}                 
\newcommand{\wkstlim}{\bw\!\bk^*\!\!-\!\!\!\lim}               
\newcommand{\Haar}[1]{\ensuremath{\sH^{_{\!a\!a\!r}}_{#1}}} 
\newcommand{\card}[1]{{\sf Card}\lb[#1\rb]}
\newcommand{\Meas}[1]{\sM \lb[#1\rb] }
\newcommand{\mtrx}[3]{\ensuremath{\lb[#1  |_{#2}^{#3} \rb]}}         
\newcommand{\Matrix}[1]{\ensuremath{\lb[\begin{array}{ccccccccccccccccccccccccr} #1 \end{array}\rb]}}   
\newcommand{\rank}[1]{\ensuremath{{\sf rank}\lb[#1\rb]}}
\newcommand{\Natur}{\dN}
\newcommand{\Zahl}{\dZ}
\newcommand{\Zahlmod}[1]{\ensuremath{\Zahl_{/{#1}}}}
\newcommand{\Cplx}{\dC}
\newcommand{\Torus}[1]{\ensuremath{{\dT}^{#1}}}  
\newcommand{\CC}[1]{\ensuremath{{\lb[ #1 \rb]}}}
\newcommand{\CO}[1]{\ensuremath{{\lb[ #1 \rb)}}}
\newcommand{\OC}[1]{\ensuremath{{\lb( #1 \rb]}}}
\newcommand{\OO}[1]{\ensuremath{{\lb( #1 \rb)}}}              
\newcommand{\boxdiagram}[8] { \begin{array}{rcl}    #1 & \stackrel{#2 }{\longrightarrow} & #3 \\    \lb. \maketall #4 \rb\downarrow & &  \lb\downarrow \maketall #5 \rb. \\    #6 & \stackrel{ #7}{\longrightarrow} &  #8 \\   \end{array} }           

\begin{document}
\ETDS{1269}{1287}{22}{2002}
\runningheads{M. Pivato, R. Yassawi}{Limit Measures for Affine
 Cellular Automata}

\title{Limit Measures for Affine Cellular Automata}
\author{Marcus Pivato\affil{1}
 and
 Reem Yassawi\affil{2}\footnote{This research partially supported by NSERC Canada.}}

\address{\affilnum{1}   Department of Mathematics, University of Houston,\\
Houston, TX 77204-3476 USA \\ \email{ {\tt pivato@math.uh.edu}
or {\tt mpivato@trentu.ca}} \\
\affilnum{2}  Department of Mathematics, Trent University,
 Lady Eaton College,\\
Peterborough, Ontario,  K9L 1Z6  Canada\\
\email{ {\tt ryassawi@trentu.ca}}}

\recd{$8$ August $2000$ and accepted in revised form $17$ September $2001$;\\
{\tt ArXiv} version submitted $31$ August $2002$}

\begin{abstract}
  Let $\dM$ be a monoid (e.g. $\Natur$, $\Zahl$, or $\Zahl^D$), and
$\sA$ an abelian group.  $\sA^\dM$ is then a compact abelian group; a
{\em linear cellular automaton} (LCA) is a continuous endomorphism
$\gF:\sA^\dM \into \sA^\dM$ that commutes with all shift maps.

   Let $\mu$ be a (possibly nonstationary) probability measure on
$\sA^\dM$; we develop sufficient conditions on $\mu$ and $\gF$ so that
the sequence $\{\gF^N \mu\}_{N=1}^\oo$ weak*-converges to the { Haar
measure} on $\sA^\dM$, { in density} (and thus, in \Cesaro average as
well).  As an application, we show: \ if $\sA=\Zahl_{/p}$ ($p$ prime),
\ $\gF$ is any ``nontrivial'' LCA on $\sA^{\lb(\Zahl^D\rb)}$, \ and
$\mu$ belongs to a broad class of measures (including most Bernoulli
measures (for $D \geq 1$) and ``fully supported'' $N$-step Markov
measures (when $D=1$), then $\gF^N \mu$ weak*-converges to Haar
measure in density.

\end{abstract}

\section{Introduction}

  Let $\sA$ be a finite set, and let $\dM$ be a monoid (e.g. $\dM =
\Zahl^D$, $\Natur^E$, or $\Zahl^D \x \Natur^E$).  Let $\sA^\dM$ be
the {\bf configuration space} of $\dM$-indexed sequences in $\sA$.
Treat $\sA$ as a discrete space; then $\sA^\dM$ is compact and totally
disconnected in the Tychonoff product topology.  The action of $\dM$ on
itself by translation induces a natural {\bf shift action} of $\dM$ on
configuration space: \ for all $\fe\in\dM$, and $\ba
\in \sA^\dM$, define $\shift{\fe}[\ba] \ = \ \mtrx{b_\fm}{\fm\in\dM}{}$
 where, $\forall \fm, \ \ \ b_\fm = a_{\fe.\fm}$, where ``$.$'' is the
monoid operator (``$+$'' for  $\dM = \Zahl^D \x \Natur^E$, etc.).

  A {\bf cellular automaton} (CA) is a continuous self-map
$\gF:\sA^\dM\into\sA^\dM$ which commutes with all shifts: for any
$\fe\in\dM$, \ \ $\gF\circ\shift{\fe} = \shift{\fe}\circ\gF$.  Hedlund
\cite{HedlundCA} proved that any such map is determined by a {\bf
local function} $\gf:\sA^\dU\into\sA$, where $\dU\subset\dM$ is some
finite set (thought of as a ``neighbourhood around the identity'' in
$\dM$), so that, for any $\ba = \mtrx{a_\fm}{\fm\in\dM}{}\in\sA^\dM$,
with $\gF(\ba) \ = \ \mtrx{b_\fm}{\fm\in\dM}{}\in\sA^\dM$, we have:
\[
\forall \fm\in\dM,  \ \ \ 
  b_\fm \ = \ \gf\lb(\ba_{|(\fm.\dU)}\rb). 
\]
  If $\sA$ is a finite abelian group with operator ``+'', then  $\sA^\dM$
is a compact abelian group under componentwise addition.  A {\bf
linear cellular automaton} (LCA) is a CA which is also a group
endomorphism from $\sA^\dM$ to itself.  This is equivalent to
requiring $\gf$ to be a group homomorphism from $\sA^\dU$ into $\sA$.
An {\bf affine} cellular automaton (ACA) is one having a local map of the
form $\gf = \gh + b$, where $\gh:\sA^\dU\into\sA$ is a homomorphism,
and $b \in \sA$ is some constant.

The term ``linear'' comes from the special case when $\sA = \Zahl_{/p}$,
 for some prime
$p$.  Since $\Zahl_{/p}$ is also a finite field, this map is actually
a {\em linear} map from the $(\Zahl_{/p})$-vector space
$(\Zahl_{/p})^\dU$ into $\Zahl_{/p}$; \ it generally takes the form:
\beqn
\gf \lb[ \ba \rb] \ = \ \sum_{\fu\in\dU} f_\fu a_\fu
\end{equation}
where $\ba = \mtrx{a_\fu}{\fu\in\dU}{}$ is an element of $\sA^\dU$,
and where $\set{f_\fu}{\fu\in\dU}$ is a set of coefficients in $\Zahl_{/p}$.

 The {\bf Haar measure} on $\sA^\dM$ is the measure $\Haar{}$ assigning mass
$A^{-N}$ to any cylinder set on $N$ coordinates, where $A=\card{\sA}$.
$\Haar{}$ is $\gF$-invariant for any LCA $\gF$, raising the question:
for what measures $\mu$ do the iterates $\gF^N \mu := \ \mu \circ
\gF^{-N}$ converge to $\Haar{}$ in the weak* topology, as $N\goto\oo$?

  This was first investigated by D. Lind \cite{Lind}, who
studied the LCA on $\lb(\Zahlmod{2}\rb)^\Zahl$ with local map
$\gf(\ba) = a_{(-1)} + a_1$.  Using methods from harmonic
analysis, Lind showed that, if $\mu$ is any nontrivial Bernoulli
 probability measure on $\lb(\Zahlmod{2}\rb)^\Zahl$, then
\[
\lim_{N \goto \oo} \frac{1}{N} \sum_{n=1}^N  \gF^{N} \mu 
\ \  = \  \ 
\Haar{}.
\]
  Lind also showed that  the sequence of
measures $\set{\gF^N \mu }{N \in \Natur}$ does not {\em itself} converge
to Haar measure; \ for all $j \in \set{2^n}{n \in \Natur}$ the measure
$\gF^j \mu$ is quite far from Haar.

   Ferrari {\em et al.}
\cite{FerrariMaassMartinez,FerMaassMartNey} studied LCA with
local maps $\gf(\ba) \ = \ k_0 a_0 + k_1 a_1$ acting on 
$\sA^\Natur$, where $\sA = \Zahlmod{q}$, \ $q = p^n$ for
some prime $p$, and $k_0$ and $k_1$ are relatively prime to $p$, and
showed \Cesaro convergence to Haar measure in the weak* topology for a
broad class of measures satisfying a certain ``correlation decay''
property, including most Bernoulli and Markov measures.  These results
are summarized in
\cite{MaassMartinezII}, where the authors also prove that most Markov
measures on $\sA^\Zahl$ will \Cesaro-converge to Haar, when
 $\sA = (\Zahl_{/2}) \dirsum (\Zahl_{/2})$, and $\gf:\sA^2 \into \sA$
is defined $\gf\lb[ (x_1,y_1), \ (x_2,y_2) \rb] = (y_1, \ x_1 + y_2)$.

  These results raise four questions: \label{four.possible.extensions}
\blist
  \item Is there some broader class of measures whose $\gF$-iterates
converge to Haar measure in \Cesaro mean?

  \item Rather than { \Cesaro} convergence, can we obtain convergence
{ in density}?  (If $\mu$ is stationary, then \Cesaro-convergence to
$\Haar{}$ is equivalent to convergence in density.  However, when
$\mu$ is nonstationary, convergence in density is a stronger result.)

  \item  For what other linear CA can we prove convergence
to Haar? What about affine CA?

  \item Can these results be generalized to LCA on higher dimensional 
lattices (e.g. $\dM=\Zahl^D$, $\Natur^D$) or 
nonabelian monoids such as free groups?
\elist

  We address these questions by developing a sufficient condition for
the sequence of measures $\seq{\gF^N \mu }{N \in \Natur}$ to converge,
{\em in density}, to Haar measure, where $\gF:\sA^\dM \into \sA^\dM$
is an LCA, and $\dM$ is a finitely generated monoid.  We require the
measure $\mu$ to have a kind of mixing property, called {\em harmonic
mixing} ---we demonstrate that, for example, Bernoulli measures (on
$\sA^\dM$, where $\dM$ is any monoid) and $N$-step Markov measures (when
$\dM$ is $\Zahl$ or $\Natur$) have this property.  We also require the
automata $\gF$ to have a kind of ``expansiveness'' property, called
{\em diffusion}, which we show is true for all ``nontrivial'' LCA
when $\sA=\Zahlmod{p}$.

\paragraph*{This paper is organized as follows:} in \S\ref{sect:prelim}, we
develop background on harmonic analysis over $\sA^\dM$
(\S\ref{sect:Harmonic}) and linear cellular automata
(\S\ref{sect:CA}).  In \S\ref{sect:Hmix} we discuss {\em harmonic mixing}
and exhibit some examples of it.  In \S\ref{sect:diffuse} we discuss
{\em diffusion} and its consequences.  In \S\ref{diffuse.lattice}, we
show that, for any prime $p$ and $D \geq 1$, if $\sA = \Zahl_{/p}$,
then all ``nontrivial'' linear cellular automata on
$\sA^{\lb(\Zahl^D\rb)}$ are diffusive; \ hence, such automata take
harmonically mixing measures on $\sA^{\lb(\Zahl^D\rb)}$ into Haar
measure.

\paragraph*{Notation:}
Elements of $\sA$ will be written as $a,b,c,\ldots$. We often
identify the elements of $\sA$ with the set $\CO{0..p} \ := \
\{0,1,\ldots,p-1\}$.
Sans-serif letters (e.g. $\fm,\fn,\fu,\ldots$) are elements of $\dM$.
Boldface letters (e.g. $\ba,\bb,\ldots$) are elements of $\sA^\dM$,
and $\ba = \mtrx{a_\fm}{\fm\in\dM}{}$.  Capitalized Gothic letters
(eg. $\gF$, $\gG$) denote cellular automata.  The corresponding
lower-case Gothic letters (eg. $\gf$, $\gothg$) denote the
corresponding local maps.

\section{Preliminaries\label{sect:prelim}}

\subsection{Harmonic Analysis on $\sA^\dM$
\label{sect:Harmonic}}

  Let $\Torus{1}$ be the unit circle group.  A {\bf character} of
 $\sA$ is a group homomorphism $\phi:\sA \into \Torus{1}$. 
Let $\hsA$ be the group of all characters of $\sA$.

If $\sA=\Zahlmod{n}$, then $\hsA$ is canonically isomorphic with $\sA$. 
First define $\bgam\in\hsA$ by
\[
\bgam(a) \ \ =  \ \ \exp\lb(\frac{2 \pi \bi}{n} a\rb).
\]
(where we identify $\sA$ with $\CO{0..n}$ in the obvious way).
Then, for each $k \in \sA$ and $a \in \sA$, define $\bgam^k\in \hsA$\ by:
$\bgam^k(a) \ := \ \bgam(k\cdot a) \ \ =  \ \ \exp\lb(\frac{2 \pi \bi}{n} k\cdot a\rb)$,
where ``$k\cdot a$'' refers to multiplication, mod $n$.
Then $\hsA = \set{\bgam^k }{k \in \sA}$, and the map
$\sA \ni k \mapsto \bgam^k \in \hsA$ is an isomorphism.

 Let $\h{\sA^\dM}$ be the group of characters of $\sA^\dM$.
If $\sA$ is any abelian group, then $\h{\sA^\dM}$
is in bijective correspondence with the set
\[
\set{ \mtrx{\chi_\fm}{\fm\in\dM}{} \ \  \in \ \  \lb(\hsA\rb)^\dM}
{ 
\chi_\fm = \chr{} \ \mbox{ for all but finitely many $\fm\in\dM$} }.
\]
If $\mtrx{\chi_\fm}{\fm\in\dM}{}$ is such a sequence,
then define 
\[ 
\bchi \ = \ \Tensor_{\fm\in\dM}\chi_\fm\  \in\ \h{\sA^\dM}.
\]
That is:\ 
if $\ba \ = \ \mtrx{a_\fm}{\fm\in\dM}{}$ is an element of $\sA^\dM$,
then $\D \bchi (\ba) \ = \ \prod_{\fm\in\dM} \chi_\fm (a_\fm)$, 
(where all but finitely many terms in this product are equal to 1.)
The sequence
 $\mtrx{\chi_\fm}{\fm\in\dM}{}$ is called the {\bf coefficient system} of
 $\bchi$.  The {\bf rank} of the character $\bchi$ is the number of
 nontrivial entries in  $\mtrx{\chi_\fm}{\fm\in\dM}{}$.

  For example, if $\sA=\Zahlmod{n}$, then $\h{\sA^\dM}$ is naturally
  isomorphic to the group
\[
\set{\maketall \mtrx{\chi_\fm}{\fm\in\dM}{} \ \  \in \ \sA^\dM }
{\chi_\fm = 0 \ \mbox{ for all but finitely many $\fm\in\dM$} }.
\]
If $\mtrx{\chi_\fm}{\fm\in\dM}{}$ is such a sequence,
then let $\D \bchi = \Tensor_{\fm\in\dM} \bgam^{\chi_\fm} :\sA^\dM \into \Torus{1}$.  Thus, if $\ba \ = \ \mtrx{a_\fm}{\fm\in\dM}{}$ is an element of $\sA^\dM$,
then
$\D\bchi (\ba) \ \ = \ \ 
\prod_{\fm\in\dM} \exp\lb(\frac{2\pi\bi}{p} \chi_\fm\cdot a_\fm\rb)$.

  Let $\Meas{\sA^\dM}$ be the space of (possibly nonstationary)
probability measures on $\sA^\dM$.  If $\mu\in\Meas{\sA^\dM}$, then
the {\bf Fourier coefficients} of $\mu$ are defined:
\[
\hmu[\chi] \ = \ \inn{\mu,\chi} \ = \ \int_{\sA^\dM} \chi \ d \mu,
\]
for all $\bchi \in \h{\sA^\dM}$.  These coefficients
completely identify $\mu$.  We will use the following basic result
from harmonic analysis:

\Theorem{\label{character.weak.star.convergence}}
{
 If $\mu_1,\mu_2,\mu_3,\ldots,\mu_\oo\in\Meas{\sA^\dM}$, then
\[
\statement{$\h{\mu_n}[\bchi] \goesto{n\goto\oo}{} \h{\mu_\oo}[\bchi]$ for all $\bchi \in
\h{\sA^\dM}$}
\iff
\statement{$\mu_n \goesto{n\goto\oo}{} \mu_\oo$ in the weak*-topology}.
\]
}

\subsection{Linear Cellular Automata \label{sect:CA}}
 
  If $\sA=\Zahlmod{n}$, and $\gf:\sA^\dU\into\sA$ is a homomorphism, then
there is a unique collection of constant
{\bf coefficients}\ $\mtrx{f_\fu}{\fu\in\dU}{}\in\sA^\dU$ so that, 
for any $\ba=\mtrx{a_\fu}{\fu\in\dU}{}\in \sA^\dU$, we have:
$\D
 \gf(\ba) \ = \ \sum_{\fu\in\dU} f_\fu a_\fu$.\ 
  Thus, if $\gF:\sA^\dM\into\sA^\dM$ is the corresponding LCA, then,
$\forall \ \ba=\mtrx{a_\fm}{\fm\in\dM}{}\in \sA^\dM$,\ \
$\D
 \gF(\ba) \ = \ \sum_{\fu\in\dU} f_\fu \cdot \shift{\fu}(\ba)$.
  In other words, we can formally write $\gF$ as a
 ``polynomial of shift maps'':
\[
 \gF \ = \ \sum_{\fu\in\dU} f_\fu \cdot \shift{\fu}.
\]
  This defines an isomorphism of between the ring of LCA over
$\sA^\dM$ and the ring of formal polynomials
with coefficients in $\sA$ and ``powers'' in $\dM$.  Composition of
cellular automata corresponds to multiplication of these polynomials.

\Proposition{\label{CA.conv}}
{ If 
$\D\gF = \ \sum_{\fu\in\dU} f_\fu \shift{\fu}$,   
and  
$\D \gG =\  \sum_{\fv\in\dV} g_\fv \shift{\fv}$,  
 then
  $\D \gF \circ \gG = \ \sum_{\fk \in \dM} 
\lb( \sum_{{\fu \in \dU, \ \fv\in\dV} \atop{\fu.\fv = \fk}} f_\fu g_\fv  \rb) \shift{\fk}$.
\qed
}

\breath

If $\gF$ is a linear cellular automata with coefficient system
$\bff \ = \ \mtrx{f_\fu}{\fu\in\dU}{}$, and $\chi$ is a character
with coefficient system $\bchi \ = \ \mtrx{\chi_\fm}{\fm\in\dM}{}$,
then define $\bchi \ast \bff \ = \  \mtrx{\xi_\fk}{\fk\in\dM}{}$,
where, for all $k \in \dM$,\ \  $\xi_k \in \sA$ is defined:
\[
\xi_k \ =  \ 
  \sum_{{\fu \in \dU, \  \fm\in\dM}\atop {\fm.\fu = \fk}} \lb( \chi_\fm \cdot f_\fu \rb)
\]
(almost all terms in this sum are equal to $0$).  Then it is not hard
to show:

\Proposition{\label{CA.char.conv}}
{
 If $\chi \in \h{\sA^\dM}$ and $\gF:\sA^\dM \into \sA^\dM$ are 
determined by coefficient systems $\bchi$ and $\bff$
respectively, then 
$\chi \circ \gF$ is also a character, and is determined by coefficient system
$\bchi \ast \bff$.
\qed}

\section{Harmonic Mixing \label{sect:Hmix}}

 A  measure $\mu$ on $\sA^\dM$ is called {\bf harmonically
mixing} if, for all $\eps > 0$, there is some $R > 0$ so that, for all
$\chi \in \h{\sA^\dM}$, \ \ $
\statement{ $\rank{\chi} > R$} \
\Implies \
   \statement{ $\lb| \hmu[\chi] \rb| < \eps$}
$.
(Notice that this definition does not require $\mu$
to be stationary.)

 For example, $\Haar{}$ is harmonically mixing;\ indeed, for all
$\chi$ except the trivial character $\chr{}$, we have: \ $\inn{\chi, \
\Haar{}} \ = \ 0$.

Let $\Meas{\sA^\dM; \Cplx}$ be the Banach algebra of complex-valued
measures (with convolution operator ``$*$'' and the total variation
norm ``$\norm{\bullet}{var}$''), and let $\sH\subset
\Meas{\sA^\dM;\Cplx}$ be the set of harmonically mixing measures.

\Proposition{\label{thm:mix.closed}}
{ $\sH$ is an ideal of $\Meas{\sA^\dM;\Cplx}$, closed under $\norm{\bullet}{var}$. }
\bthmprf
  $\sH$ is clearly closed under linear operations.  To show that $\sH$
is a convolution ideal, use the fact that $\h{\mu\ast\nu} =
\hmu \cdot \hnu$ and that $\hnu$ is bounded by $\norm{\nu}{var}$.  Thus, if $\mu$ is harmonically mixing, then so are $\mu\ast\nu$ and $\nu\ast\mu$.

 To show closure in $\norm{\bullet}{var}$, use the fact that,
for any measures $\mu$ and $\nu$, $ \norm{\mu-\nu}{var} \ = \ \sup
\set{\lb|\inn{\phi,\mu}-\inn{\phi,\nu}\rb|}
{\phi\in\bC\lb(\sA^\dM;\Cplx\rb), \ \ \norm{\phi}{\oo} = 1}$.
\ethmprf

  Not all measures on $\sA^\dM$ are harmonically mixing.  For example,
 $\fm \in \dM$, we say $\mu\in\Meas{\sA^\dM}$ is {\bf
 $\fm$-quasiperiodic} if there is an orthonormal basis $\seq{\bxi_n}{n
 \in \Natur}$ of $\bL^2\lb(\sA^\dM; \ \mu\rb)$, consisting entirely of
 eigenfunctions of the shift map $\shift{\fm}$.  It is not difficult
 to show that, if $\mu$ is $\fm$-quasiperiodic for any $\fm
 \not=\Id{\dM}$, then $\mu$ is not harmonically mixing.

\breath

 Also, if $\bchi:\sA^\dM \into \Torus{1}$ is a nontrivial character, then the
{\bf Markov subgroup} \cite{KlausSchmidt,LindMarcus,KitS}
induced by $\bchi$ is defined:
\[
\sA^\dM_{\bchi} \ := \ \set{\ba \in \sA^\dM}{\bchi \circ \shift{\fm}(\ba) \ = \ 1, \ \forall \fm\in\dM}.
\]
If $\sA^\dM_{\bchi}$ is nontrivial, it is a subshift of finite type.
If $\mu$ is a stationary probability measure on $\sA^\dM_\chi$, then
$\mu$ cannot be harmonically mixing: \ if $\fm_1,\ldots,\fm_K \in \dM$
are spaced widely enough apart, and $\D\bxi \ := \ \prod_{k=1}^K \bchi
\circ \shift{\fm_k}$, then $\inn{\mu,\ \bxi} \ = \ 1$, no matter how
large $K$ becomes.

 However, $\mu$ may still be harmonically mixing relative to the
elements of $\h{\sA^\dM_{\bchi}}$; \ see Corollary \ref{subgroup.haar}.

\subsection{Harmonic Mixing of Bernoulli Measures \label{bernoulli.mix}}

\Proposition{}
{
  Let $\sA = \Zahl_{/p}$, where $p$ is prime.  Let $\beta$ be any
measure on $\sA$ which is {\em not} entirely concentrated on one
point.  Let
$\D  \beta^{\tensor\dM} \ = \ \Tensor_{\fm\in\dM} \beta$ 
be the corresponding Bernoulli measure on $\sA^\dM$.
Then $\beta^{\tensor\dM}$ is harmonically mixing.
}
\bthmprf
 $\forall k \in \sA$, let $c_k := \inn{\bgam^k, \ \beta}$,
where $\bgam^k\in\hsA$ is as in \S\ref{sect:Harmonic}.
Since $p$ is prime,
$|c_k| < 1$, unless $k= 0$, while $c_0 = 1$. 
Thus, $\displaystyle c := \max_{0< k < p} |c_k| <1$.  
Thus, if $\bchi \in \h{\sA^\dM}$ and $\rank{\bchi} = R$, then
\ $\D
\lb|\inn{\bchi, \ \beta^{\tensor\dM} }\rb|
\ = \
\lb|\inn{\Tensor_{\fm\in\dM} \chi_\fm, \ \Tensor_{\fm\in\dM} \beta}\rb| 
\ = \
\lb|\prod_{\fm\in\dM} \inn{ \chi_\fm, \beta} \rb|
\ < \ c^R$ \
becomes arbitrarily small as $R$ gets large.
\ethmprf

  A similar argument shows:

\Proposition{\label{bernoulli.mix.2}}
{
  Let $\sA$ be an arbitrary finite abelian group, and
$\beta$ a measure on $\sA$.  Suppose that, for any subgroup $\sG \subset \sA$,
the support of $\mu$ extends over {\em more than one} coset of $\sG$.
Then $\beta^{\tensor\dM}$ is harmonically mixing.
\qed}

\Corollary{}
{
 $\sH$ is weak* dense in $\Meas{\sA^\dM}$.
}
\bthmprf
Let $\mu\in\Meas{\sA^\dM}$ be arbitrary.  For any $\eps\in\CC{0,1}$,
let $\nu_\eps= \beta_\eps^{\tensor \dM} \in\Meas{\sA^\dM}$ be the
Bernoulli measure with one-dimensional marginal $\beta_\eps$, where
$\beta_\eps[0] = 1-\eps$, and, for all $a\in\sA\setminus\{0\}$,
$\beta_\eps[a] =
\eps/(A-1)$  (where $A=\card{\sA}$).   $\nu_\eps\in\sH$ by
Proposition \ref{bernoulli.mix} so
$\nu_\eps
\ast\mu\in\sH$ also, by Proposition \ref{thm:mix.closed}.

  We want to show that $\D \wkstlim_{\eps\goto0} \nu_\eps \ast\mu \ = \ \mu$;
\ it is equivalent to show that  $\D\lim_{\eps\goto0} \h{\nu_\eps \ast\mu}
 \ = \ \hmu$,\ pointwise.   Clearly, $\D \wkstlim_{\eps\goto0} \nu_\eps \ = \ \del_0$, where $\del_0$ is the point mass on the constant
zero configuration $0\in\sA^\dM$.   Thus, for any $\bchi\in\h{\sA^\dM}$,
\ $\D \h{\nu_\eps\ast\mu}(\bchi)
\ \ = \ \
\h{\nu_\eps}(\bchi)\cdot\hmu(\bchi)
\ \ \goesto{\eps\goto0}{} \ \
\h{\del_0}(\bchi)\cdot \hmu(\bchi)
\ \ =\  \ \bchi(0)\cdot\hmu(\bchi)\ \  = \ \ \hmu(\bchi)$.
\ethmprf

\subsection{Harmonic Mixing of Markov Measures}

  Now let $\dM = \Zahl$, and suppose that $\mu$ is a {\bf Markov
measure} on $\sA^\dZ$, is determined by the transition probability
matrix $\brbQ = \mtrx{q^a_b}{a,b \in \sA}{}$, with stationary
probability vector $\bnu \ = \ \mtrx{\nu_a}{a \in \sA}{}$, so that
$\brbQ \cdot \bnu \ = \ \bnu$.   Thus, if $\bc \in \sA^\Zahl$
is a $\mu$-random configuration, then \ $\bnu_a \ = \  \Prob{c_0 \ = \ a }$,
 \ and  \ $q^a_b \ :=   \ \Prob{c_1 = b \ | \ c_0 = a}$.

\Proposition{\label{thm:mix.markov}}
{Let $\sA$ be any finite abelian group.  If all entries of 
$\brbQ$ are nonzero, then
$\mu$ is harmonically mixing.
}
\bthmprf Let $\Cplx^\sA$ be the set of all functions $\xi:\sA \into \Cplx$.

  Define the operator $\sQ:\Cplx^\sA \into \Cplx^\sA$ as follows:
for any $\xi \in \Cplx^\sA$ and any $a\in\sA$,
\[
\sQ[\xi](a) \ = \ \sum_{b \in \sA} q^a_b \xi(b).
\]
In other words, \ \ $\sQ[\xi](a) \ = \ \inn{ \xi, \bq^a }$,
where $\bq^a$ is the ``$a$th'' column of the matrix $\brbQ$, and
we treat $\xi$ as an $\sA$-indexed vector.

\breath

  Next, for any $\chi \in \hsA$, define the multiplication-by-$\chi$
operator: $\sM_\chi:\Cplx^\sA \into \Cplx^\sA$ so that, for
any $\xi \in \Cplx^\sA$ and any $a\in\sA$, \ \ \ 
$\sM_\chi[\xi](a) \ = \ \chi(a) \cdot \xi(a)$. 

Now, suppose $\bchi = \chi_0 \tensor \chi_1 \tensor \ldots \tensor
 \chi_N$ is a character on $\sA^\dM$ (in other words, $\D\bchi =
 \Tensor_{n\in\Zahl}\chi_\fn$, but $\chi_\fn = \chr{}$ for all $n > N$ and
 $n < 0$).

\Claim{ 
\setcounter{enumi}{1}
\begin{list}{(\alph{enumii})}{\usecounter{enumii}}
\item  \label{inn.chi.mu} If $N \geq 1$, then 
$\ \  \D  \inn{\bchi,\mu} 
\ = \ \inn{ \sM_{\chi_0} \circ \sQ \circ 
\sM_{\chi_1} \circ \sQ \circ \ldots \circ \sM_{\chi_{N-1}} \circ \sQ [\chi_N]
, \ \ \bnu}.$
\item \label{norm.mult}
 For any $\phi \in \Cplx^\sA$ and any $\bchi \in \hsA$,  \ \ 
$\norm{\sM_\chi[\phi]}{\oo} \ = \ \norm{\phi}{\oo}$.
\item \label{norm.Q}
 For any {\em non}constant $\phi \in \Cplx^\sA$, \ \ \ 
$\norm{\sQ[\phi]}{\oo} \ < \ \norm{\phi}{\oo}.$ 
\end{list}
}
\bclaimprf (a) is just linear algebra.  (b) is because, $\forall a\in\sA$,
\ $\lb|\chi(a)\rb|=1$.  To see (c), note that,
 $\forall a \in \sA$,
\ $\D
\lb|\sQ[\phi](a)\rb|
 \  = \ 
   \lb| \sum_{b \in \sA} q_b^a \phi(b)  \rb| 
 \ \leq  \  
    \sum_{b \in \sA} q_b^a \lb| \phi(b) \rb| 
 \ \leq \ \
    \sup_{b \in \sA} \lb| \phi(b) \rb| \ \ = \ \ \norm{\phi}{\oo}$. \
The first (triangle) inequality is an equality if and only if all
the elements of $\set{\phi(b)}{b \in \sA}$ have the same phase angle.
The second inequality is an equality if and only if they all have the
same magnitude.  Hence, $ \lb|\sQ[\phi](a)\rb| \leq \norm{\phi}{\oo}$,
with equality if and only if $\phi$ is constant. \eclaimprf

 Now, for any $\xi,\zeta \in \hsA$, with $\zeta \not=\chr{}$, define
$\sP_{\xi,\zeta} := \sM_{\xi} \circ \sQ \circ \sM_{\zeta} \circ \sQ$.
Then $\sP_{\xi,\zeta}:\Cplx^\sA \into \Cplx^\sA$ is a linear operator.
If  $\Cplx^\sA$ is endowed with the $\norm{\bullet}{\oo}$ norm,
then let $\norm{\sP_{\xi,\zeta}}{\oo}$ be the operator norm of
$\sP_{\xi,\zeta}$.

\Claim{\label{norm.P}  $\norm{\sP_{\xi,\zeta}}{\oo} \ < \ 1$.}
\bclaimprf
  Let $\phi \in \Cplx^\sA$, with $\norm{\phi}{\oo} = 1$.
If $\phi$ is {\em not} constant, then by {\bf Claim}
\ref{norm.Q},
$\norm{\sQ[\phi]}{\oo} < 1$; \ thus, by {\bf Claim} \ref{norm.mult} and
{\bf Claim} \ref{norm.Q}, \ \   $ \norm{\sM_{\xi} \circ \sQ \circ \sM_{\zeta} \circ \sQ[\phi]}{\oo} \leq \norm{\sQ[\phi]}{\oo} < 1$.
  If $\phi$ {\em is} constant,
 then $\sM_{\zeta}\circ\sQ[\phi]$ is {\em not} constant;\  thus, by
{\bf Claim} \ref{norm.Q}, \ \ 
$\norm{\sM_{\xi} \circ \sQ \circ \sM_{\zeta} \circ \sQ[\phi]}{\oo} 
 \ < \norm{ \sM_{\zeta} \circ \sQ[\phi]}{\oo} \leq \norm{\phi}{\oo} = 1$.

  $\Cplx^\sA$ is finite-dimensional, so the unit ball $\dB$ 
relative to the supremum norm $\norm{\bullet}{\oo}$ is compact;  hence
$\D
\norm{\sP_{\xi,\zeta}}{\oo}\ = \ 
   \sup_{\phi \in \dB}  \norm{\sP_{\xi,\zeta}[\phi]}{\oo} \ < \ 1.
$
\eclaimprf
  Thus, for all $\xi,\zeta \in \hsA$, with $\zeta \not=\chr{}$, let
$c_{\xi,\zeta} := \norm{\sP_{\xi,\zeta}}{\oo}$,
and let
\[
C := \max \set{c_{\xi,\zeta}}
  {{\xi,\zeta \in \hsA}\  \mbox{and} \ {\zeta\not=\chr{}}}
\]
Thus, since $c_{\xi,\zeta} < 1$ for all $\xi,\zeta$, and since
$\hsA$ is finite, we conclude that $C < 1$ also.  So,
given any $\eps > 0$, if $K$ is large enough, then $C^K < \eps$.

  Now, if $\rank{\bchi} > 2K$, then the product:
\ \ \  $\bchi \ = \ \chi_0 \tensor \chi_1 \tensor \ldots \tensor \chi_N$
can be rewritten:
\beq
\bchi & = &
 \lb(\underbrace{\chr{}\tensor\ldots\tensor\chr{}}_{n_0}\rb)
\tensor \lb(\xi_1 \tensor \zeta_1\rb) \tensor
\lb(\underbrace{\chr{}\tensor\ldots\tensor\chr{}}_{n_1}\rb)
\tensor \lb(\xi_2 \tensor \zeta_2\rb) \tensor
\ldots \\
&& \ldots  \tensor
\lb(\underbrace{\chr{}\tensor\ldots\tensor\chr{}}_{n_{R-1}}\rb)
\tensor \lb(\xi_R \tensor \zeta_R\rb) \tensor
\lb(\underbrace{\chr{}\tensor\ldots\tensor\chr{}}_{n_R}\rb), 
\eeq
where $R> K$, and, for all $r \in \CO{0..R}$,\  \ $\xi_r,\zeta_r$ are successive elements in the list
$\chi_0,\chi_1,\ldots,\chi_{N-1}$, with $\zeta_r \not=\chr{}$,
and where $n_0,n_1,\ldots,n_R \geq 0$, so that
\ \ \ $n_0 + n_1 + \ldots + n_R + 2R \ = \ N$. \ \ 
Thus, the operator $\sM_{\chi_0} \circ \sQ \circ 
\sM_{\chi_1} \circ \sQ \circ \ldots \circ \sM_{\chi_{N-1}} \circ \sQ$
\ can be rewritten as \newline  $\lb(\sQ^{n_0} \circ \sP_{\xi_1,\zeta_1}\rb)
\circ  \lb(\sQ^{n_1} \circ \sP_{\xi_2,\zeta_2}\rb)  \circ \ldots
\circ  \lb(\sQ^{n_{R-1}} \circ \sP_{\xi_R,\zeta_R} \rb) \circ \sQ^{n_R}$.
But then
\beq
\lb|\inn{\bchi, \mu}\rb|
& =_{(1)} &
 \lb|\inn{ \sM_{\chi_0} \circ \sQ \circ 
\sM_{\chi_1} \circ \sQ \circ \ldots \circ \sM_{\chi_{N-1}} \circ \sQ[\chi_N],
 \ \ \bnu}\rb| \\
& \leq_{(2)} &
\norm{  \sM_{\chi_0} \circ \sQ \circ 
\sM_{\chi_1} \circ \sQ \circ \ldots \circ \sM_{\chi_{N-1}} \circ \sQ[\chi_N]}{\oo} \\
& \leq_{(3)} &
\norm{  \sM_{\chi_0} \circ \sQ \circ 
\sM_{\chi_1} \circ \sQ \circ \ldots \circ \sM_{\chi_{N-1}} \circ \sQ}{\oo} \\
& = &
\lb \| \lb(\sQ^{n_0}  \circ \sP_{\xi_1,\zeta_1}\rb)
\circ  \lb(\sQ^{n_1} \circ \sP_{\xi_2,\zeta_2}\rb)  \circ \ldots 
\rb. \\
&&\hspace{9em}
\lb. 
 \ldots \circ  \lb(\sQ^{n_{R-1}} \circ \sP_{\xi_R,\zeta_R} \rb) \circ \sQ^{n_R} \rb\|_\oo \\
& \leq &
\norm{ \sP_{\xi_1,\zeta_1}}{\oo} \cdot
\norm{\sP_{\xi_2,\zeta_2}}{\oo} \cdot \ldots \cdot
\norm{\sP_{\xi_R,\zeta_R}}{\oo}\\
& \leq &
  C^R \ \ < \ \ C^K \ \ <  \ \ \eps 
\eeq
(1)  by {\bf Claim} \ref{inn.chi.mu}.\ \   (2)  $\bnu$ is a probability measure.
\ \ (3)  $\norm{\chi_N}{\oo}  = 1$.

 In summary, if $\rank{\bchi} > 2R$, then  \ \  
$\lb|\inn{\bchi, \mu}\rb| < \eps$.
\ethmprf

\Corollary{\label{thm:mix.markov.abscont}}
{Let $\sA$ and $\mu$ be as in Theorem \ref{thm:mix.markov}, and
suppose $\nu$ is a measure on $\sA^\Zahl$ absolutely continuous
relative to $\mu$.  Then $\nu$ is also harmonically mixing.}
\bthmprf
  Let $\D\phi = \frac{d\nu}{d\mu}$, and suppose first that $\D\phi =
\frac{\chr{[\ba]}}{\mu[\ba]}$ is the (renormalized) characteristic
function of some cylinder set $[\ba] = \set{\bb\in\sA^\Zahl}{\bb_{\dU}
\ = \ \ba}$, where $\dU = \CC{-U\ldots U}\subset \Zahl$ and $\ba\in\sA^\dU$.
Thus $\nu = \mu_{[\ba]}$, the (renormalized) restriction of $\mu$
to a probability measure on $[\ba]$ (that is: $\nu(B) =
\mu\lb([\ba]\intsct B\rb)/\mu\lb([\ba]\rb)$ for any measurable
$B\subset\sA^\Zahl$).

   Let $\D \bchi = \Tensor_{n=-N}^N \chi_n$ be a character,
and suppose $N>U$.  Let $\D \bchi_{(-)} = \Tensor_{n=-N}^{-U-1} \chi_n$
 and $\D \bchi_{(+)} =  \Tensor_{n=U+1}^{N} \chi_n$.
Let $\mu_{[\ba]}^{(+)}\in\Meas{\sA^\OO{U...\oo}}$ be the projection of $\mu_{[\ba]}$ onto coordinates $\OO{U...\oo}$  (thus,
if $\bb \in \sA^\OC{U..N}$, then \ $\D
\mu_{[\ba]}^{(+)}[\bb] \ = \ 
q^{a_U}_{b_{(U+1)}} \cdot q^{b_{(U+1)}}_{b_{(U+2)}} \cdot \ldots \cdot
q^{b_{(N-1)}}_{b_{N}}$).
Similarly, let  $\mu_{[\ba]}^{(-)}$ be the projection of
$\mu_{[\ba]}$ onto coordinates $\OO{-\oo\ldots -\!U}$.
Thus, using the Markov property of $\mu$,
\[
\inn{\bchi,\nu}
\ = \
\inn{\bchi_{(-)},\ \mu_{[\ba]}^{(-)}} \cdot
\lb(\prod_{u=-U}^U \chi_u {a_u}\rb)
\cdot
\inn{\bchi_{(+)},\ \mu_{[\ba]}^{(+)}}
\]
Now, analogous to Claim \ref{inn.chi.mu} of Theorem \ref{thm:mix.markov}, we have:
\[
\inn{\bchi_{(-)},\mu_{[\ba]}^{(-)}} 
\ = \ \inn{ \sM_{\chi_{(-N)}} \circ \sQ \circ 
\sM_{\chi_{(1-N)}} \circ \sQ \circ \ldots \circ \sM_{\chi_{(-U-1)}} \circ \sQ [\chi_U]
, \ \ \bq_{a_{(-U)}}  },
\]
where $\bq_{a_{(-U)}}$ is the $a_{(-U)}$th ``row'' of transition matrix
$\brbQ$, and, in a manner analagous to the proof of Theorem
\ref{thm:mix.markov}, we can show that
\[
\lb|\inn{\bchi_{(-)},\ \mu_{[\ba]}^{(-)}} \rb| \goto 0 \ \ 
\mbox{ as} \ \ 
\rank{\bchi_{(-)}}\goto\oo.
\]  By a similar argument (with reversed
time), we can show
\[
\lb|\inn{\bchi_{(+)},\ \mu_{[\ba]}^{(+)}} \rb| \goto
0 \ \ \mbox{as} \ \ \rank{\bchi_{(+)}}\goto\oo.
\]
This shows that $\nu$ is harmonically mixing.

  The case when $\phi$ is {\bf simple} ---ie. a finite linear
combination of characteristic functions of cylinder sets ---then
follows immediately, via Proposition \ref{thm:mix.closed}.  If $\phi
\in \bL^1(\mu)$ is arbitrary, let $\seq{\phi_n}{n\in\Natur}$ be a
sequence of simple functions converging to $\phi$ in the $\bL^1$ norm.
Let $\seq{\nu_n}{n\in\Natur}$ be the corresponding measures
(all harmonically mixing); \  thus, $\seq{\nu_n}{n\in\Natur}$
converges to $\nu$ in total variation norm, so $\nu$ is also
harmonically mixing, by Proposition \ref{thm:mix.closed}.
\ethmprf

  Notice that the measure $\nu$ need not be stationary (and will not be,
unless $\phi$ is shift-invariant.)

  An {\bf $N$-step Markov process} is analogous to a Markov process,
but the probability distribution of each letter is dependent upon the
previous $N$ letters, and conditionally independent of what comes
before.  When $N=0$, we have a Bernoulli process; \ when $N=1$, a
standard Markov process.  In general, an $N$-step process is
determined by a collection of transition probabilities $\brbQ \ = \ 
\set{q^\ba_b}{\ba\in\sA^\CO{0..N}, \ \ b \in \sA}$ so that, for each
$\ba\in\sA^\CO{0..N}$, \ $\D \sum_{b\in\sA} q^\ba_b = 1$.  One can then
find a (generally unique) stationary distribution $\bnu$ on
$\sA^\CO{0..N}$.  
\Corollary{\label{thm:mix.N-step.markov}}
{ 
Let $\sA$ be a finite abelian group and let $\alp$ be an $N$-step
Markov process on $\sA$, where all elements of $\brbQ$ are nonzero.
Then $\alp$ is harmonically mixing.
}
\bthmprf
  Let $\sB = \sA^\CC{1..N}$, and consider the standard {\em $N$-block coding}
map  $\phi:\sA^\Zahl\into\sB^\Zahl$, defined:
\[
\phi(\ldots,a_1,\ldots,a_N,a_{N+1},\ldots,a_{2N},\ldots)
\ = \ \lb( \ldots, \Matrix{a_1\\\vdots\\a_N}, \ 
\Matrix{a_{N+1}\\ \vdots\\ a_{2N}}, \ldots \rb)
\]
This is an isomorphism of topological groups, and the following diagram
commutes:
\[
\boxdiagram{\sA^\Zahl}{\shift{N}}{\sA^\Zahl}
	{\phi} 			{\phi}
	{\sB^\Zahl}{\shift{}}{\sB^\Zahl}
\]
Thus, $\beta = \phi^*\alp$ is a (1-step) Markov measure,
with transition matrix $\brbP = \mtrx{p^\ba_\bb}{\ba,\bb \in \sB}{}$,
where $p^{(a_1,\ldots,a_N)}_{(b_1,\ldots,b_N)} = 
q^{(a_1,\ldots,a_N)}_{b_1} \cdot q^{(a_2,\ldots,a_N,b_1)}_{b_2} 
\cdot q^{(a_3,\ldots,a_N,b_1,b_2)}_{b_3} \cdot \ldots \cdot
q^{(a_N,b_1,\ldots,b_{(N-1)})}_{b_N}$.  Clearly, if all entries
of $\brbQ$ are nonzero, then, so are all entries of
$\brbP$, and thus, by Proposition \ref{thm:mix.markov}, $\beta$ is
harmonically mixing.  Hence, it suffices to show:

\claim{If $\beta$ is harmonically mixing, then so is $\alp$.}
\bclaimprf
 The isomorphism $\phi:\sA^\Zahl\into\sB^\Zahl$ induces
isomorphism $\hphi:\h{\sB^\Zahl}\into\h{\sA^\Zahl}$ given:
$\hphi(\bchi):=\bchi\circ\phi$.  Thus, $\hphi^{-1}:\h{\sA^\Zahl}\into\h{\sB^\Zahl}$ is an isomorphism, and, for any $\bchi\in\h{\sA^\Zahl}$,
\blist
 \item $\D\rank{\hphi^{-1}(\bchi)} \geq \frac{1}{N}\rank{\bchi}$.
 \item $\inn{\hphi^{-1}(\bchi), \ \beta} \ = \ \inn{\bchi,\ \alp}$.
\elist
Thus, if $\rank{\bchi}$ is large, then so is $\rank{\hphi^{-1}(\bchi)}$;
then $\inn{\hphi^{-1}(\bchi), \ \beta}$ is small, and thus so is
$ \inn{\bchi,\ \alp}$.
\eclaimprf
\ethmprf

   An $N$-step Markov measure satisfying the hypothesis of Corollary
\ref{thm:mix.N-step.markov} gives nonzero probability to all
cylinder sets of finite length; we might say it has ``full support''.

   The technique of Proposition \ref{thm:mix.markov} can be used to prove
the corresponding result for stationary Markov random fields on free
groups and monoids \cite{Spitzer,Zachary}.  Now, instead of {\em one} transition
probability matrix, there are several: one for each generator of the
group/monoid.  As long as there are finitely many generators, the
bound $C$ on the operator norms in the proof of Proposition
\ref{thm:mix.markov} will still be less than 1, and the same argument
can be used to show:

\Theorem{\label{thm:mix.markov.field}}
{
 Let $\sA$ be a finite group.  Let $\dM$ be a free group or free
  monoid on finitely many generators, and let $\mu$ be a stationary
  Markov random field on $\sA^\dM$ so that all entries in all
  transition probability matrices are nonzero. 
  Then $\mu$ is harmonically mixing.
\qed}

\section{Diffusive Linear Automata \label{sect:diffuse}}

  Let $\gF:\sA^\dM \into \sA^\dM$ be a linear cellular automaton.
We say that $\gF$ is {\bf diffusive} if, for every nontrivial
$\chi \in \h{\sA^\dM}$, \ \ $\D \lim_{n \goto \oo}
 \rank{ \chi \circ \gF^n }\ \ = \ \  \oo$.

For example, let $\sA = \Zahl_{/p}$ for some prime $p$.  Let $\dM$ be
 the free group or free monoid on $D\geq2$ generators.  Let $\dU
 \subset \dM$ be a finite subset with at least two elements, so that
 each $\fu\in\dU$ is a product of at least two distinct generators.
It is straightforward to show:
\[
\mbox{\sl If $0 \not= f_\fu \in \Zahl_{/p}$, for all $\fu\in\dU$,
then the automaton} \ \
 \gF \ = \ \sum_{\fu\in\dU} f_\fu \shift{\fu} 
\ \ 
  \mbox{ is diffusive.}
\]
   Unfortunately, linear cellular automata on $\lb(\Zahlmod{p}\rb)^{\lb(\Zahl^D\rb)}$ are never diffusive: \ if $\gF$ is such an LCA (e.g. $\gF = 
f_0 + f_1\cdot \shift{\fm_1} + f_2\cdot \shift{\fm_2}$), then,
for any $n\in\Natur$
the LCA $\gF^{\lb(p^n\rb)}$ is simply $\gF$, ``rescaled'' by a factor of $p^n$
(e.g. $\gF^p = 
f_0 + f_1\cdot \shift{(p\cdot \fm_1)} + f_2\cdot \shift{(p\cdot \fm_2)}$).
This follows from the {\bf Fermat property} for the field $\Zahlmod{p}$.
Thus, the polynomial of $\gF^{\lb(p^n\rb)}$ has the same number of nonzero
coefficients as that of $\gF$, so $\gF$ cannot ``diffuse'' along the subsequence
$\set{p^n}{n\in\Natur}$.

 This motivates a slight weakening of the
concept of diffusion:  we say that $\gF$ is {\bf diffusive in density} if,
 for every nontrivial $\chi \in \h{\sA^\dM}$, there is a subset
$\dJ \subset \Natur$ of \Cesaro density 1 so that $\D
\lim_{{j \goto \oo} \atop {j \in \dJ}} \rank{ \chi \circ \gF^j }
 \ \ = \ \ \oo$.

\Theorem{\label{diffuse.and.mix.converge}}
{
 Let $\sA$ be a finite abelian group, and $\dM$ a countable monoid.
 Suppose that $\gF:\sA^\dM \into \sA^\dM$ is a linear cellular
 automata, and suppose that $\mu$ is a
 measure on $\sA^\dM$ that is harmonically mixing.
\blist
\item If $\gF$ is diffusive, then \ 
$\D \wkstlim_{j \goto \oo} \gF^j \mu
\ \ =  \  \Haar{}$. 

\item If $\gF$ is diffusive {\em in density}, then there is
 a  set $\dJ \subset \Natur$ of \Cesaro density 1
so that
 \   $\D \wkstlim_{{j \goto \oo} \atop {j \in \dJ}} \gF^j \mu
\  = \  \Haar{}$. \  
Thus\ 
$\D \wkstlim_{N\goto\oo} \frac{1}{N} \sum_{n=1}^N \gF^n \mu
\  =  \  \Haar{}$.
\elist
}
\bthmprf We'll prove convergence in density, from which \Cesaro convergence
follows immediately.  The proof of strict convergence is much the same.
 
   We'll show that the Fourier coefficients of
$\gF^n \mu$ all converge to zero in density.  Weak* convergence in
density then follows by Theorem \ref{character.weak.star.convergence}.
So, let $\chi \in \h{\sA^\dM}$. Then
\[
\inn{\chi, \ \gF^n \mu}
\ \  = \ \ 
\int_{\sA^\dM} \chi \  d \lb( \gF^n \mu\rb)
 \ \  = \ \ 
\int_{\sA^\dM} \chi\circ \gF^n \  d \mu 
\ \  = \ \  
 \inn{\lb( \chi\circ \gF^n\rb) , \mu }
\]
Now, since $\gF$ is diffusive in density, we can find a subset
$\dJ_\chi \subset \Natur$ of density 1 so that\ 
$\D
\lim_{{j \goto \oo} \atop {j \in \dJ_{\chi}}}  \rank{ \chi \circ \gF^j }
\ \ = \ \  \oo$.
But then, since $\mu$ is harmonically mixing, it follows that\ 
$\D
\lim_{{j \goto \oo} \atop {j \in \dJ_{\chi}}} \inn{\lb( \chi\circ \gF^j\rb) , \mu }
 \ \ = \ \ 0$.

Now, $\sA$ is finite and $\dM$ is countable; \ thus,
 $\h{\sA^\dM}$ is countable, so we can find a ``common tail set''
$\dJ \subset \Natur$ so that:
\bitem
  \item $\dJ$ has \Cesaro density 1.
  \item For every $\chi \in \h{\sA^\dM}$, there is some $N > 0$ so that
$\dJ_\chi \intsct \CO{N..\oo} \ \subset \ \dJ$.
\eitem
  (\cite{Petersen}, Remark 2.6.3, or 
\cite{KakutaniJones}).  Thus, for all  $\chi \in \h{\sA^\dM}$,\
$\D\lim_{{j \goto \oo} \atop {j \in \dJ}}  
\inn{\lb( \chi\circ \gF^n\rb) , \mu }
\ \ = \ \ 0$.
\ethmprf

  The same reasoning applies stationary measures supported on
shift-invariant subgroups of $\sA^\dM$:

\Corollary{\label{subgroup.haar}}
{
 Let $\gF:\sA^\dM \into \sA^\dM$ be as in Theorem 
\ref{diffuse.and.mix.converge}.  Suppose that $\sG \subset \sA^\dM$ is a
closed subgroup, and $\mu$ is a measure supported on $\sG$.  Suppose
$\mu$ is harmonically mixing relative to the elements of $\h{\sG}$, and
$\gF(\sG) \subseteq \sG$.   If $\gF$ is diffusive (in density), then\
$\D \wkstlim_{j \goto \oo} \gF^j \mu
\ \ =  \  \Haar{\sG}$,\  
 where $\Haar{\sG}$ is the Haar measure on the compact group $\sG$
(and where convergence is either absolute or in density, as appropriate).
\qed}

\breath

  To extend these results to {\em affine} cellular automata,
use the following:

\proposition{\label{converge.on.subgroup}}
{ 
Let $\sA$ be any finite group.  Let $\gF:\sA^\dM \into \sA^\dM$ be
 an LCA with local transformation $\gf:\sA^\dU \into \sA$.  Let $c \in
 \sA$ be some constant, and let $\gG$ be the ACA with local
 map $\gothg:\sA^\dU \into \sA$ given: $\gothg(\ba) \ = \
 \gf(\ba) + c$.

  Let $\mu$ be a measure on $\sA^\dM$, and let $\dJ \subset \Natur$. 
\[
\mbox{If} \ \ 
\wkstlim_{{j \goto \oo} \atop {j \in \dJ}} \gF^j \mu 
\ \ = \ \ \Haar{} 
  \ \ 
\mbox{ then } \ \ 
 \wkstlim_{{j \goto \oo} \atop {j \in \dJ}} 
  \gG^j \mu
\ \ = \ \  \Haar{}. 
\]
}
\bthmprf
 Let $\bc_0 \in \sA^\dM$ be the constant configuration  whose
entries are all equal to $c$, and, $\forall n\in\Natur$,
let $\bc_n = \gF^n(\bc_0)$.  Let $\Bh_n = \bc_0+\bc_1+\ldots+\bc_n$,
and define $\gH_n:\sA^\dM\into\sA^\dM$
by: $\gH_n(\ba) = \ba + \Bh_n$.  A simple computation establishes:
\[
\forall n\in\Natur,  \ \ \ \gG^n \ = \ \gH_n\circ\gF^n.
\]
 If $\bchi \in \h{\sA^\dM}$, then for any
$\ba \in \sA^\dM$, we have:
$\bchi \circ \gH_n(\ba) \ = \ \bchi\lb(\ba + \Bh_n\rb) \ = \ 
K_n \cdot \bchi(\ba)$,
where $K_n = \bchi\lb(\Bh_n\rb)$ is some element of $\Torus{1}$.  Concisely:  $\bchi \circ \gH_n \ = \ K_n \cdot \bchi$.  

  Now, $\gF^n \mu$ converges in density to the Haar measure, in the
weak* topology, which is equivalent to saying: for every nontrivial
character $\bchi$, there is a subset $\dJ \subset \Natur$ of
density one such that $\D \lim_{{j \goto \oo} \atop {j \in \dJ}}  
\inn{ \bchi \circ \gF^j, \ \mu } \ = \ 0$.

Now, for any $j$, \ \ $
\inn{ \bchi \circ \gG^j, \ \mu }
\ \  = \ \  
 \inn{ \bchi \circ  \gH_j \circ \gF^j, \ \mu } 
\ \  = \ \ 
\inn{  K_j \cdot \bchi \circ \gF^j,\  \mu} 
\ \  = \ \ 
 K_j \cdot \inn{ \bchi \circ \gF^j,\  \mu}$. \ 
 But $\lb|K_j\rb| = 1$, and thus,
\[
\lb|\inn{ \bchi \circ \gG^j, \ \mu }\rb| \ \ = \ \ 
\lb| \inn{ \bchi \circ \gF^j,\  \mu } \rb| \ \ \goesto{j \goto \oo}{j \in \dJ}
\ \ 0.\]
Since this is true for each character, we conclude that
$\gG^n \mu$ also converges in density to the Haar measure.
\ethmprf

\section{Diffusion on Lattices \label{diffuse.lattice}}

  Say that an LCA $\gF$ on $\sA^\dM$ is {\bf
nontrivial} if $\gF$ is {\em not} merely a shift map or the identity map.
If we write $\gF$ as a polynomial of shift maps, then $\gF$ is 
nontrivial if this polynomial contains two or more nonzero
coefficients.

\Theorem{ \label{unboundedthm}  }
{ Let $p$ be a prime number, and $\sA = \Zahl_{/p}$.  Let $D \geq 1$.  Then
  any nontrivial linear cellular automaton on $\sA^{\lb(\Zahl^D\rb)}$
  is diffusive in density.  }

\breath

  The proof of this theorem will occupy the rest of this section.
We will eventually accomplish a reduction to the case when $D=1$; \
hence, the reader may initially find it helpful to assume
$D=1$, and to treat all elements of $\Zahl^D$ (indicated as
vectors, eg. ``$\vm$'') as elements of $\Zahl$ instead (indicated
as scalars, eg. ``$m$'').  It will also be helpful to first work through
the details of the proof in the special case when $p=2$; \  we will
make reference to this special case in footnotes.

  We will represent LCA using the polynomial notation introduced
in \S\ref{sect:CA}.  It will be convenient to write these polynomials
in a special recursive fashion.  For example, suppose $D=1$, and
suppose that $g_0,g_1,g_2
\in \CO{1..p}$, and $\ell_0,
\ell_1, \ell_2 \in \Zahl$.  Let $\sG$ be the linear CA on $\sA^\Zahl$
defined: \  $\sG = g_0
\shift{\ell_0} \ + \ g_1 \shift{\ell_1} \ + \ g_2 \shift{\ell_2}$.
Then we can rewrite $\gG$ as:
\beqn
\label{J.is.2}
 \gG  \ = \ g_0 \cdot \lb( \gF \circ \shift{\ell_0} \rb), \ \ \mbox{ where} \ \
\gF \ \ = \ \ \Id{} \ + \ f_1 \shift{m_1} \lb(1 \ + \ f_2 \shift{m_2}\rb) \, ,
\end{equation}
with $m_1 = \ \ell_1 - \ell_0$, \ $m_2 = \ell_2 - \ell_1$, \
$f_1 = g_0^{-1} g_1$ and  $f_2 = g_1^{-1} g_2$ (with inversion in
the field $\Zahl_{/p}$)\footnote{If $p=2$, we can assume that $f_1 = f_2 = 1$.}.
More generally, we have the following:

\Lemma{}
{  Let $g_0,g_1,\ldots, g_J \in \CO{1..p}$, and $\vec{\ell}_0, \vec{\ell}_1, \ldots, \vec{\ell}_J \in \Zahl^D$, and suppose
that $\gG$ is the linear CA on $\sA^{(\Zahl^D)}$ defined:
\beqn
\gG \ = \ g_0 \shift{\vec{\ell}_0} + g_1 \shift{\vec{\ell}_1} + \ldots + g_J \shift{\vec{\ell}_J}\, .
 \label{automaton}
\end{equation}
Then $ \gG  \ = \ g_0 \cdot \lb( \gF \circ \shift{\vec{\ell}_0} \rb),$  \ 
where:
\begin{eqnarray}
\label{multipolynomial}
 \lefteqn{\gF  \ = \ } \\
 && \Id{} \ + \ f_1 \shift{\vm_1} \lb( \Id{} \ + \ f_2 \shift{\vm_2}
   \lb[ \maketall \ldots \lb( \Id{} + 
 f_{J-1} \shift{\vm_{J-1}} \lb[ \Id{} \ + \ f_J \shift{\vm_J}\rb] \rb) \ldots  \rb] \rb)  \nonumber
\end{eqnarray}
 and, for all $j \in \CC{1..J}$, \ \ $\vm_j \ = \ 
 \vec{\ell}_j  -\vec{\ell}_{j-1}$, and $ f_j \  = \ g_{j-1}^{-1} \cdot g_j$.
\qed}

\breath

   Composing with the shift $\shift{\vec{\ell}_0}$ and multiplying by the
scalar $g_0$ does not affect the diffusion property; \
hence, it
is sufficient to prove Theorem for polynomials like (\ref{multipolynomial}).
On first reading, it may be helpful to assume that $J=2$, as in
(\ref{J.is.2}).

  By Proposition \ref{CA.conv}, the powers $\gF^N$ of the linear
cellular automaton $\gF$ correspond to powers of the corresponding
polynomial.  To prove Theorem \ref{unboundedthm}, we will therefore
need to develop some machinery concerning multiplication of
polynomials with coefficients in $\Zahl_{/p}$, by using the classical
formula of Lucas \cite{Lucas} for the mod $p$ binomial coefficients,
which has become ubiquitous in the theory of LCA.

\Definition{$p$-ary expansion, Index set.}
{

If $n \, \in  \Natur $, 
then the {\bf $p$-ary expansion} of $n$\ is the sequence\  $\dP(n) = 
\{ n^{[i]} \}_{i=0}^{\infty} \ \in \CO{0..p}^\Natur$, such that
$\D n = \sum_{i=0}^{\infty} n^{[i]} p^{i}$.
 
The {\bf index set} $\sS(n)$ of $n$ is defined to be: \ \ \
$\sS(n) = \set{ i \, \in \dN}{ n^{[i]} \neq  0 }.$
}

 If $m\in\Natur$, then let $[m]_p$ be the congruence class of $m$, mod $p$.

\paragraph*{Lucas' Theorem:}
{\sl
 Let $N,n \in \Natur$, with $p$-ary expansions as before.  Then
\[
\lb[ N \atop n \rb]_p \ = \ \prod_{k=0}^\oo \lb[ N^{[k]} \atop n^{[k]} \rb]_p
\]
where we define
  $\D\lb(0 \atop 0 \rb) \  = \  1$, \ 
and \ $\D\lb(a \atop b \rb) \ \  = \ \  0$, for any $b > a > 0$.
}
\qed\cite{Lucas}
\breath

 Write ``$n \ll  N$'' if $n^{[i]} \leq N^{[i]}$ for all $i \in \Natur$.
Thus, Lucas' Theorem implies:
\[
\lb( \ \lb[ N \atop n \rb]_p \not= \ 0 \ \rb) \ 
\iff \statement{ $n \ll  N$ }.
\]

If $n \, \in \Natur $, then the {\bf Lucas set} of $n$ is the set\ \ 
$\sL (n) \ :=\  \set{ k  \in \Natur}{k \ll  n  }$.

  The following elementary arithmetic observation will be used later.

\Lemma{\label{add.binary}}
{ 
Let $n_{1}$, $n_{2}, \ldots , n_{L} \in \Natur$.  If $M>0$, and
for all $\ell \in \CC{1..L}$ and $i \geq M $, $n_{l}^{[i]} = 0$, then,
for all $i \geq M+ \lceil \log_{p} L \rceil$, \ \ 
 $\D\left( \sum_{l=1}^{L} n_{l} \right )^{[i]} \ \ = \ \ 0$.\ 

(Here, $\lceil n \rceil \, \in \dN$ is the smallest integer such that $\lceil n \rceil \geq n$.)
}
\bthmprf
The ``$\log_p$'' term comes from the fact that, in summing $L$
distinct $p$-ary numbers, there is the possibility of up to
$\log_p[L]$ digits of carried value spilling forward.
\ethmprf
 
  With Lucas' theorem, one can obtain expressions for powers of linear
automata.  For example, if $f \in  \Zahl_{/p}$, and $\gF$ is the linear
automata on $\sA^\Zahl$ defined by $\gF(x) = x + f\cdot \shift{} (x)$, then
Lucas' Theorem tell us that
\[
\gF^{N} (x) = \sum_{k \, \in \sL ( N)} 
  \pbinom{N}{k}{p} f^k\cdot \shift{k} x \, .
\]
Next, if $m_1,m_2 \in \Zahl$, and $f_1,f_2 \in \Zahl_{/p}$, and
$\gF$ is as in (\ref{J.is.2}), then
\beq
\gF^{N} & = &
\sum_{ k_1 \, \in \sL (N) }
 \pbinom{N}{k_1}{p}  f_1^{k_1}\cdot \shift{m_1 k_1} \lb(1 + f_2\cdot \shift{m_2}\rb)^{k_1} \\
& = & 
 \sum_{ k_1 \, \in \sL (N) } \pbinom{N}{k_1}{p} f_1^{k_1}\cdot \shift{m_1 k_1} \ 
\lb ( \sum_{{k_2} \,\in \sL (k_1)} \pbinom{k_1}{k_2}{p} f_2^{k_2} \cdot \shift{m_2 k_2} \rb ) \\
& = &
   \sum_{k_1 \in \sL(N)} \sum_{ {k_2} \in \sL(k_1)} 
    \lb( \pbinom{N}{k_1}{p} \pbinom{k_1}{k_2}{p} f_1^{k_1} f_2^{k_2}\rb) \cdot  \shift{m_1 k_1 +m_2 k_2} .\\ 
& = &
   \sum_{k_1 \in \sL(N)} \sum_{ {k_2} \in \sL(k_1)} 
    f_{(k_1,k_2)} \cdot  \shift{m_1 k_1 +m_2 k_2} ,\\ 
&&\mbox{  where we define} \ \  f_{(k_1,k_2)} \ \  := \ \  \pbinom{N}{k_1}{p} \pbinom{k_1}{k_2}{p} f_1^{k_1} f_2^{k_2}.
\eeq
A similar argument works in $\Zahl^D$, and for an arbitrary number
of ``nested'' polynomial terms of this type.  This leads to the
following

\Lemma{\label{power.of.nested.poly}}
{  If $\vm_1,\vm_2,\ldots,
\vm_J \in \Zahl^D$, and $f_1,\ldots,f_J \in \Zahl_{/p}$, and $\gF$
is as in (\ref{multipolynomial}), then

\vspace{-1em}

\begin{eqnarray}
\gF^N
  & = &  
\sum_{\bk \in \sL^J(N)} f_{(\bk)} \shift{\inn{\bk,\Bm}}, \label{multipoly.exp}\\
\mbox{where:} \hspace{4em}
 \Bm
  & := & \lb[\vm_1,\vm_2,\ldots,\vm_J\rb], \nonumber \\
  \sL^J(N)
  & := & \set{  \lb[k_1,k_2,\ldots,k_J\rb] \in \Natur^J}
		{ k_J \ll k_{J-1} \ll \ldots k_2 \ll k_1 \ll N }, \nonumber
\end{eqnarray}

 and, for any such $\bk \ = \ \lb[k_1,k_2,\ldots,k_J\rb]$, we define
\beq
\inn{\bk,\Bm} & := & k_1 \vm_1 + k_2 \vm_2 + \ldots + k_J \vm_J \\
\mbox{and} \ \
  f_{(\bk)} & := &  
\pbinom{N}{k_1}{p} \pbinom{k_1}{k_2}{p} \ldots   \pbinom{k_{J-1}}{k_J}{p} 
f_1^{k_1} f_2^{k_2} \ldots f_J^{k_J}.\\
\eeq
(the dependence on $N$ is suppressed in the notation ``$f_{(\bk)}$''.)
\qed}

\bprfof{Theorem \ref{unboundedthm}}
 It suffices to prove the
theorem for polynomials $\gF$ like (\ref{multipolynomial}).
So, suppose $\gF$ is {\em not} diffusive in density.  Thus, there exists
some nontrivial character $\bchi \in \h{\sA^{(\Zahl^D)}}$, some
$R \, \in \Natur$, and a subset $\dB \subset \Natur$
(of ``bad'' numbers), of upper density $\delta > 0$, so that, for all $n \, \in \dB$, \ \ \  $\rank{\bchi \circ \gF^{n}} \leq R$.   

\breath

  Now, for each $\vn \in \Zahl^D$, let $\pr{\vn}:\sA^{(\Zahl^D)} \into
  \sA$ be projection onto the $\vn$th coordinate: \ \ $\pr{\vn} (\ba)
  \ = \ a_\vn$.  Let $\bgam:\sA\into\Torus{1}$ be the character
introduced in \S\ref{sect:Harmonic}: \ 
$\bgam(a) \ = \ \exp\lb(\frac{2 \pi\bi}{p} \cdot a\rb)$.  Thus,
there is a finite subset $\sQ \subset \Zahl^D$, and a collection of
coefficients $\set{\chi_\vq \in \CO{1..p}}{\vq \in \sQ}$ so that
$\bchi$ is defined\footnote{In the case when $p=2$, we can write this:\ \
$\D\bchi (\ba) = \prod_{\vq \in \sQ} (-1)^{a_{\vq}}$.}:
\beqn
\bchi (\ba) = \prod_{\vq \in \sQ}  \bgam\lb( \chi_{\vq} \cdot \pr{\vq}(\ba)\rb)
\label{character}
\end{equation}
Thus, if $\gF^N$ is as in (\ref{multipoly.exp}) of Lemma
\ref{power.of.nested.poly}, and $\bchi$ is as in (\ref{character}),
then, by Proposition
\ref{CA.char.conv},  the character $\bchi
\circ \gF^{N}$ has the following expansion\footnote{When $J=2=p$,
and $D=1$
the expansion is:
\[
\bchi \circ \gF^{n} (\bx)
\ \  = \ \ \prod_{q \in \sQ} \
\prod_{k_1 \in \sL(n)} \ \prod_{ k_2 \in \sL(k_1)} \
(-1)^{(k_1 m_1+k_2 m_2 +q)}
\].}:
\beqn
\bchi \circ \gF^{N} 
\ \  = \ \ \prod_{\vq \in \sQ} \
\prod_{\bk \in \sL^J(N)}  \
 \bgam\lb( \chi_\vq \cdot f_{(\bk)} \cdot \pr{\lb(\inn{\bk,\Bm} +\vq\rb)}\rb) .
\label{expansion}
\end{equation}
Note that, for every $\vq \in \sQ$ and $\bk \in \sL^J(N)$, the
factor $\bgam\lb(\chi_\vq \cdot f_{(\bk)} \cdot \pr{\lb(\inn{\bk,\Bm} +\vq\rb)}\rb)$
is nontrivial: \ $f_{(\bk)}$ is never a multiple of $p$, and thus, if
$\chi_{\vq}$ is nontrivial, then $\bgam \lb(\chi_\vq \cdot f_{(\bk)} \cdot \pr{\lb(\inn{\bk,\Bm} +\vq\rb)}\rb)$
is also nontrivial.  Thus, the
{\em only} way the coefficients of the character defined by
(\ref{expansion}) can be trivial is if two terms of the form
$\bgam\lb(\chi_\vq^* \cdot f_{(\bk^*)}\cdot\pr{(\inn{\bk^*,\Bm} +\vq^*)}\rb)$ and
$\bgam\lb(\chi_\vq \cdot f_{(\bk)}\cdot \pr{(\inn{\bk,\Bm} +\vq)}\rb)$ cancel
out, which can only occur when
\beqn
\label{multicancelation}
 \inn{\bk^*,\Bm} +\vq^* \ \ = \ \ \inn{\bk,\Bm} +\vq.
\end{equation}
This is an equation of $D$-tuples of integers, and hence, is only 
true if, for all $d \in \CC{1...D}$,
\beqn
\label{cancelation}
 \inn{\bk^*,\Bm}_{(d)} +q_{(d)}^*  \ \ = \ \ \inn{\bk,\Bm}_{(d)} +q_{(d)}
\end{equation}
where the subscript ``$(d)$'' refers to the $d$th component of the
$D$-tuple.

\breath

\subparagraph*{The idea of the proof} is thus as follows: In order for the rank
of the character $\bchi \circ \gF^{N}$ (for $N \in \dB$) to be less than
$R$, most of the terms in the expression (\ref{expansion}) must
cancel out; \ this requires a specific kind of ``destructive
interference'' between the the index sets $\sS(N)$ and various
translations of $\sS(N)$ so that virtually all elements $(\bk,\vq) \in \sL^J(N)
\x \sQ$ must be paired up as in equation (\ref{multicancelation}), so as to
cancel with each other. 

   Our goal, then, is to show that the equation
(\ref{multicancelation}) is hard to achieve, so that, after the dust
settles, more than $R$ nontrivial coefficients remain.  We will show
that, the set $\dB$ (indeed, {\em any} set of nonzero density) must
contain numbers for which sufficient cancellation fails to occur.

\subparagraph*{Reduction to Case $D=1$:}
 In order for cancellation of terms (ie. equation
(\ref{multicancelation})) to occur in $\Zahl^D$, equation
(\ref{cancelation}) must be true for {\em every} $d \in \CC{1..D}$
simultaneously.  Hence, it is enough to disrupt the equation in one
dimension.  Hence, at this point, we can reduce the argument to the
case when $D=1$.  We will treat $m_1,\ldots,m_J$ as elements of
$\Zahl$, and $\Bm = [m_1,\ldots,m_J]$ as a $J$-tuple of integers; \
thus, for any other $J$-tuple $\bk = [k_1,\ldots,k_J]$, we have
$\inn{\bk,\Bm} \ = \ k_1m_1 + \ldots k_Jm_J$.
Likewise, $\sQ$ will be some finite subset of $\Zahl$.  

\subparagraph*{Gaps in the Index set:} We will use an ergodic 
argument to show that any subset of $\Natur$ of nonzero density must
contain numbers $N$ possessing large ``gaps'' in their index sets: \ i.e.
$\dP[N]$ has long blocks of $0$'s terminated by $1$'s.  We can then
find elements $k_1^*
\in \sL(N)$ also exhibiting these long gaps.  The gap in such a
$k_1^*$ is long enough that it is impossible to find some other element
$(\bk,q) \in \sL^J(N) \x \sQ$ so that the terms in
the expression $\inn{\bk,\Bm} +\vq$ sum together to ``cancel'' the
terminating $1$ in the gap of $\sS(k_1^*)$.
 
   Since there are many of these gaps, there are many such
elements $k_1^*$, and thus, there will be at least $(R+1)$ distinct
$1$'s that remain uncancelled, and thus at least $(R+1)$ nontrivial
terms in expression (\ref{expansion}), contradicting the hypothesis
that $\rank{\bchi \circ \gF^{N}} \leq R$ for all $N \in \dB$.

\breath

 We can assume that, when we transformed
expression (\ref{automaton}) into expression (\ref{multipolynomial}),
we had $\ell_0 < \ell_1 < \ldots < \ell_J$;\ 
hence, we can assume that $m_1,\ldots,m_J \ > \ 0$.   Thus,
they have well-defined Lucas sets, $\sS(m_1),\ldots,\sS(m_J).$
So, to begin, define:
\[
\Gamma \ \ := \ 
  \max\lb[\Union_{j=1}^J \sS(m_j)\rb] 
\ + \ 
     \lb\lceil \log_p\lb( \sum_{j=1}^J  \card{\sS(m_j)} \rb) + \log_p(J)
 \rb\rceil
 +2
\]
$\Gamma$ stands for ``gap'', and is the size of the gaps we will
require. 

  Let $q_1$ be the smallest element of $\sQ$, and define
\[
\sQ_1 \ := \ \set{ \maketall q - q_1}{q \in \sQ}, \ \ 
\mbox{and} \ \ 
\sU \ := \ \Union_{q \in \sQ_1} \sS\lb(q\rb).
\]

Next, let $\bw$ be the element of \  $\CO{0..p}^{\Gamma+1}$ 
defined:  $\D\bw \ := \ (\underbrace{0,\ldots,0}_\Gamma, \ 1 )$.
We will be concerned with the frequency of occurrence of $\bw$ in
the $p$-ary expansions of integers.

\subparagraph*{Notation:}
 If  $\bs \ =\ s_{0}s_{1}\ldots s_{n}$ is a string in $\CO{0..p}^n$, then 
we define the {\bf frequency} of the word $\bw$ in $\bs$, denoted by $fr[\bw,\bs]$, as
\[
fr[\bw,\bs] 
:= 
\frac{\card{ \{ i \, \in [0,n-1]: s_{i}s_{(i+1)}....s_{(i+\Gamma)} = \bw \} }}
{n}
\]

If $ s_{i}s_{(i+1)} \ldots s_{(i+\Gamma)} = \bw$, we'll say that $\bw$ {\bf occurs} at $s_{i}$.

\claim{
\label{erglemma}
 For any $\epsilon > 0$, there
exists $M^{*}$ such that, for any $M > M^{*}$, there is a set
$\sG^M_\bw(\epsilon) \subset \CO{1...p}^M$ so that:
\[
\card{\sG^M_\bw(\epsilon)}  \ > \  (1 - \epsilon){p}^{M},
\ \ \mbox{and}, \ \ \forall \bg \in \sG^M_\bw(\epsilon),
\ \ \ fr[\bw,\bg] > \frac{(1-\epsilon )}{p^{(\Gamma+1)}}.
\]
}
\bclaimprf  Consider the ergodic dynamical system $\lb(\CO{0..p}^\Natur,\
\Haar{}, \ \shift{}\rb)$, where $\Haar{}$ is the Haar measure and
$\shift{}: \CO{0..p}^{\dN} \rightarrow \CO{0..p}^{\dN}$ is the shift action.
The set $\set{\ba \in \CO{0..p}^\Zahl}{\ba_{\CC{0...\Gamma}} = \bw}$ has measure
$p^{-\Gamma-1}$.  The result now follows from Birkhoff's Ergodic Theorem.
\eclaimprf

In particular, let
$\D \epsilon \ := \ \frac{\delta}{2p}.
\ \ \mbox{Also, let}  \ \epsilon^{*}  \ := \ \frac{1- \epsilon}{2} \Haar{}(\bw) \ = \ \frac{1- \epsilon}{2} p^{-1-\Gamma}$.

\claim{\label{M.and.N} There exist $M$ and $N$
such that the following conditions are satisfied:
\begin{enumerate}
\item
$M \epsilon^{*} > R +2,$
\item
$\sU \subset [0,M{\epsilon}^{*}]$,
\item
$N \, \in \dB \intsct [0,p^{M})$,
\item 
$fr[\bw, \vN] \ > \ (1 - \epsilon) \, p^{-1-\Gamma}$.
\end{enumerate}
}
\bclaimprf  $\dB$ has upper density $\delta$, so there is
some sequence  $\{n_{k}\}_{k=0}^\oo$ such that,
\[
\frac{\card{\dB \intsct [0,n_{k}] }}{n_{k}} \ \goesto{k\goto\oo}{} \ 
 \delta.
\]
 Find $K$ so that, for $k > K$, \ 
 $\D\frac{\card{\dB \intsct \CC{0..n_k}}}{n_k} \ > \ \frac{\del}{2}$.
Then choose $M$ large enough to satisfy [1] and [2], and such that 
$p^{M-1} \leq n_{k} \leq p^{M}$.  Thus,
\begin{eqnarray}
\card{\dB \intsct \CC{0, p^{M}}} 
& \geq & \card{\dB \intsct \CC{0, n_{k}}}
\ \  \geq\ \  \frac{\delta n_{k}}{2}
\ \  \geq \ \ \frac{ \delta p^{M-1}}{2} \nonumber \\
& \geq & \frac{\delta p^{M} }{2p}.  \label{B.bound}
\end{eqnarray}
 Also, by Claim \ref{erglemma}, let $M$ be large enough so that
there is a subset $\sG^M_\bw(\epsilon) \subset \CO{0..p}^{M}$
so that
\begin{eqnarray}
 \card{\sG^M_\bw(\epsilon)} & > & 
  (1- \epsilon) p^{M} \ \ = \ \ \lb(1- \frac{\delta}{2p}\rb) p^{M},
\label{G.bound} \\
 \mbox{and} \ \ 
 fr[\bw, \ba]  & > & (1 - \epsilon) \, p^{-1-\Gamma}, \ \ \mbox{for all $\ba \in \sG^M_\bw(\epsilon)$.} \nonumber 
\end{eqnarray}
  Now, if $\dG := \set{n \in \CC{1..p^M}}{\dP(n) \in \sG^M_\bw(\epsilon)}$,
then $\card{\dG} = \card{\sG^M_\bw(\epsilon)}$.  Thus,
combining (\ref{B.bound}) and (\ref{G.bound}), we see that
$\card{\dB \intsct \CC{0, p^{M}}} + \card{\dG} \ > \ p^M$; \ \ hence,
the two sets must intersect nontrivially.  Let $N \in 
\dB \intsct \CC{0, p^{M}} \intsct \dG$; \ then N satisfies [3] and [4].
\eclaimprf

\claim{  Let $Q = \lceil M\epsilon^{*} \rceil$.  Then
$\bw$ occurs more than $R$ times in the string: 
$ (N^{[Q+1]}\, N^{[Q+2]}\,  N^{[Q+3]}\ldots  N^{[M]} )$.}
\bclaimprf
$N$ satisfies condition [4] of Claim \ref{M.and.N}, and of course $\bw$ occurs at most 
$Q$ times in the string 
$(N^{0}\, N^{[1]}\,  \ldots N^{[Q]} ).$
Thus, beyond position $Q$, $\bw$ must occurs at least
\beq
\lb((1-\epsilon) p^{-1-\Gamma} M\rb) \  -\  Q & \geq & 
 \lb( (1-\epsilon) p^{-1-\Gamma} M\rb)  \  - \ M{\epsilon}^{*} - 1 \\
& = & (1-\epsilon) p^{-1-\Gamma} M  \ - \ 
		\lb( M \frac{1-\epsilon}{2} p^{-1-\Gamma} \rb) - 1 \\
& = & M \frac{1-\epsilon}{2} p^{-1-\Gamma} - 1
\ \ = \ \ M{\epsilon}^{*} - 1 \\
\eeq
times and so, by condition [1] of Claim 1, at least $R +1$ times.
\eclaimprf

  Say $\bw$ occurs at some positions $N^{[j_{1}]},\ N^{[j_{2}]},\
\ldots, N^{[j_{R+1}]}$ beyond $Q$.  Thus, for each  $r \in \CC{1...R+1}$,
we have: $N^{[j_r + k]} = 0$ for $0 \leq k < \Gamma$ and   
$N^{[j_r + \Gamma]}  = 1$.   In particular,
\beqn
\label{big.gaps}
 \forall r \in \CC{1...R +1}, \ \ \ 
 p^{j_r + \Gamma} \, \in  \sL (N).
\end{equation}
Now, $N \in \dB$, so $\rank{\bchi \circ \gF^{N}} \, \leq
R$. This means that in the expression (\ref{expansion}), all but at
most $R$ of the terms are cancelled by a like term. In other words,
for all but $R$ of the elements: $(\bk^*,q^*) \in \sL^J(N) \x \sQ$,
there exists some $(\bk,q) \in \sL^J(N) \x \sQ$ so that
\beqn
\label{nihil}
\inn{\bk^*,\ \Bm} \ +\ q^{*} \ \  =  \ \   \inn{\bk,\ \Bm} \ +\ q  .
\end{equation}
---we say that  $(\bk^*,q^*)$ is {\em annihilated} by $(\bk,q)$.

  However, there are $R+1$ elements in
the set $\{j_r\}_{r=1}^{R+1}$, and thus, there are $R+1$ pairs of the form
$(\bk_r^*,q_1)$, where $\bk_r^* \ = \ \lb(p^{(j_r + \Gamma)}, 0 , \ldots, 0 \rb)$.   Hence
there exists some $r$ such that the pair  $(\bk_r^{*},q_1)$
is annihilated by some other pair $(\bk,q)$.   Define
$n \ \ := \ \ j_r + \Gamma$; \ then $ \inn{\bk_r^*,\ \Bm} \ = \ m_1p^{n}$,
so we can rewrite (\ref{nihil}) as:
\beqn
\label{annihilate}
m_1p^{n} = \inn{\bk,\Bm} \ \  + \ \ (q - q_1) ,
\end{equation}
where $\bk \ = \ \lb[k_1,\ldots,k_J\rb]$ is some other
element in $\sL^J(N)$.

\claim{\label{k.bound} For all $j \in \CC{1..J}$, and all $i \geq n - \Gamma$, \
we have:
  \ ${k_j}^{[i]} \ = \ 0$.}
\bclaimprf
   First we'll show $k_1^{[i]} = 0$ for $i\geq n$.
The RHS and LHS of  (\ref{annihilate}) must come from
different terms of the expansion (\ref{expansion}), which means that
either $q \not= q_1$ or $k_j \not=0$ for some $j>1$ ; \ either way, one of the {\em
other} terms on the RHS is positive besides ``$m_1 k_1$'', and
therefore,  $m_1 k_1\ < \ m_1 p^{n}$.  Thus, 
$k_1\ < \ p^{n}$, and thus, $k_1^{[i]} = 0$ for all $i\geq n$. 

  Next we'll show $k_1^{[i]}$ for $n-\Gamma\leq i< n$.  Recall that $k_1 \, \in \sL
 (N)$, and by hypothesis, $N^{[i]} = 0$ for all $i \in
\CO{j_r...(j_r+\Gamma)}$, where $n = j_r + \Gamma$ and $n-\Gamma =
j_r$.  Thus, ${k_1}^{[i]} = 0$ for $n - \Gamma \leq i < n$.
 
  Since $k_J \ll k_{J-1} \ll \ldots \ll k_2 \ll k_1 $, the same holds for
$k_2,\ldots,k_J$.
\eclaimprf

\claim{\label{km.bound} For all $j \in \CC{1..J}$, and all
 $i \ \geq\  n-2-\log_p(J) $,\ \ we have:
 $ (m_j k_j)^{[i]} \ = \ 0$.}
\bclaimprf  Fix $j \in \CC{1..J}$.  
For any $s \in \sS(m_j)$, it follows from Claim \ref{k.bound} that
\[
\lb(m_j^{[s]} p^s k_j\rb)^{[i]} \ = \ 0, \ \ \ 
 \mbox{for all $i \geq n - \Gamma + s$.}
\]
Hence, by Lemma \ref{add.binary},
\beq
\lb( m_j k_j\rb)^{[i]} 
&  = & 
\lb( \sum_{s \in \sS(m_j)} m_j^{[s]} p^s k_j \rb)^{[i]} \\
 & = & 0,
   \ \ \forall i \ \geq
   \ n - \Gamma + \max\lb[\sS(m_j)\rb] + \log_p\lb(\card{\sS(m_j)}\rb).\\
\eeq
The claim now follows from the definition of $\Gamma$.
\eclaimprf

\claim{\label{q.bound} For all $i \geq n-3$, \ \ $(q - q_1)^{[i]} = 0 $. }
\bclaimprf
By definition, $n-3  \ > \ n-\Gamma \ = \ j_r \ > \ Q \ \geq \ M\epsilon^*$. 
Recall that condition [2] defining $M$ was: \ 
$\sU \subset [0,M{\epsilon^*}]$.  Thus,
$(q - q_1)^{[i]} = 0 $ for $i \geq M{\epsilon^*}$.
\eclaimprf

\claim{\label{final.bound} For all $i \geq n - 1$, \ \ 
$\left( \inn{\bk,\Bm} \ \ +\ \ ( q - q_1) \right)^{[i]} \ \ = \ \ 0 $.}
\bclaimprf
  Recall that $\inn{\bk,\Bm} \ = \ (k_1 m_1) + \ldots  + (k_J m_J)$; \
thus, it follows from Claim \ref{km.bound} and Lemma \ref{add.binary} that
\ \ $\inn{\bk,\Bm}^{[i]} \ \ = \ \ 0, \ \ \forall i \ \geq \ n-2.$

Thus, the claim follows from Claim \ref{q.bound} and Lemma \ref{add.binary}.
\eclaimprf

  Now, by hypothesis, $\inn{\bk,\Bm} \ + \ ( q - q_1) \ = \ m_1 p^n$.
Hence, $\lb(\inn{\bk,\Bm} \ + \ ( q - q_1)\rb)^{[i]} $  $ = \ (m_1 p^n)^{[i]}$ for
all $i \in \Natur$.  In particular, if $I \ := \ \min\lb[\sS(m_1)\rb]\geq 0 $, then
\[
\lb(\inn{\bk,\Bm} \ + \ ( q - q_1)\rb)^{[I+n]} \ = \ 
 (m_1 p^n)^{[I+n]} \ = \ m_1^{[I]} \ \not= \ 0.
\]
But $I + n  \ \geq \  n$, so this is a contradiction of Claim \ref{final.bound}.
{\tt \hrulefill $\Box$}
\eprf

\section*{Conclusion}

  We have shown that harmonically mixing measures on $\sA^\dM$, when
acted upon by diffusive linear cellular automata (possibly with an
affine part), will weak*-converge, in \Cesaro mean, to the Haar
measure.  This sufficient condition is broadly applicable: \ in
particular, if $\sA = \Zahl_{/p}$, then any nontrivial linear cellular
automata acting upon a ``fully supported'' $N$-step Markov measure (in
$\Zahl$, a regular tree, or a free group) or a nontrivial Bernoulli
measure (in $\Zahl^D$) will converge to Haar in \Cesaro mean.

    In a forthcoming paper
\cite{PivatoYassawi2}, we generalize the results on diffusion to the case when
$\sA=\Zahlmod{n}$ ($n\in\Natur$ arbitrary) and
$\sA=\lb(\Zahlmod{(p^r)}\rb)^J$ ($p$ prime, $r,J\in\Natur$), and we
demonstrate harmonic mixing for Markov random fields on
$\sA^{\lb(\Zahl^D\rb)}$, for $D\geq2$.  However, many questions remain
unanswered.  What other classes of measures on $\Zahl$ or $\Zahl^D$
are harmonically mixing?  Measures on $\sA^\dM$ exhibiting
quasiperiodicity cannot be harmonically mixing; \ what is the \Cesaro
limit of such a measure, if anything?  Also, what LCA are diffusive,
when $\dM$ is neither a lattice nor a free group?

\paragraph*{Acknowledgments:}

  We would like to thank David Poole of Trent University for
introducing us to Lucas' Theorem, and Dan Rudolph of the University of
Maryland for reminding us that, for stationary $\mu$,
\Cesaro convergence to $\Haar{}$ is equivalent to convergence in density. 

\bibliographystyle{plain}
\bibliography{../LaTeX/bibliography}

\end{document}